\documentclass{article}
\usepackage{amssymb}
\usepackage{mathrsfs}
\usepackage[all]{xy}
\usepackage{ifthen}
\usepackage[english,francais]{babel}
\usepackage[leftbars]{changebar}
\usepackage{vmargin}

\setpapersize{A4}
\setmarginsrb{37mm}{20mm}{30mm}{37mm}{10mm}{6mm}{0mm}{10mm}

\newcommand{\tiret}{\rule[0.6ex]{1.3ex}{0.26ex}}

\newcommand{\spacebeforeenv}{\vspace{1ex}}
\newcommand{\spaceafterenv}{\vspace{1ex}}

\newenvironment{traitsurlecote}
{\cbstart
\setcounter{changebargrey}{0}      
}
{\cbend
}

\newcommand{\mysquare}{\rule[0.3ex]{1ex}{1ex}}
\newcommand{\myqed}{\hfill \mysquare}
\newcommand{\myqedenv}{}

\newcommand{\envfont}{\sf\bfseries}

\newcounter{ctheo}
\renewcommand{\thectheo}{\arabic{ctheo}}

\newcommand{\metaenvironnementhm}[2]  
{
  \refstepcounter{ctheo}
  \ifthenelse{ \equal{#2}{toto} }{
    \spacebeforeenv \begin{traitsurlecote} \noindent {\envfont #1 \thectheo}
  }{
    \spacebeforeenv \begin{traitsurlecote} \noindent {\envfont #1 \thectheo{} \tiret\, #2.}
  }
}

\newenvironment{proposition}[1][toto]
{
  \metaenvironnementhm{Proposition}{#1}
}
{\end{traitsurlecote}\spaceafterenv}

\newenvironment{lemma}[1][toto]
{ \metaenvironnementhm{Lemma}{#1} }
{\end{traitsurlecote}\spaceafterenv}

\newenvironment{lemme}[1][toto]
{ \metaenvironnementhm{Lemme}{#1} }
{\end{traitsurlecote}\spaceafterenv}

\newenvironment{corollaire}[1][toto]
{ \metaenvironnementhm{Corollaire}{#1} }
{\end{traitsurlecote}\spaceafterenv}

\newenvironment{definition}[1][toto]
{ \metaenvironnementhm{Definition}{#1} }
{\end{traitsurlecote}\spaceafterenv}

\newenvironment{keyword}
{\noindent  {\envfont Keywords.~}}
{\spaceafterenv}

\newenvironment{msc}
{\noindent  {\envfont Mathematics Subject Classification (2000).~}}
{\spaceafterenv}

\newcommand{\metaenvironnement}[2]  
{
  \refstepcounter{ctheo}
  \ifthenelse{ \equal{#2}{toto} }{
    \spacebeforeenv \noindent {\envfont #1 \thectheo}
  }{
    \spacebeforeenv \noindent {\envfont  #1 \thectheo{} \tiret\, #2.}
  }
}

\newcommand{\metaenvironnementsimple}[2]  
{
  \ifthenelse{ \equal{#2}{toto} }{
    \spacebeforeenv \noindent {\envfont #1}      
  }{
    \spacebeforeenv \noindent {\envfont #1 \tiret\, #2.}
  }
}

\newenvironment{example}[1][toto]
{ \metaenvironnement{Example}{#1} }
{\myqedenv\spaceafterenv}

\newenvironment{exemple}[1][toto]
{ \metaenvironnement{Exemple}{#1} }
{\myqedenv\spaceafterenv}

\newenvironment{remarque}[1][toto]
{ \metaenvironnement{Remarque}{#1} }
{\myqedenv\spaceafterenv}

\def\og{\leavevmode\raise.3ex\hbox{$\scriptscriptstyle\langle\!\langle$~}}
\def\fg{\leavevmode\raise.3ex\hbox{~$\!\scriptscriptstyle\,\rangle\!\rangle$}}

\newcommand{\numtoi}[1]{
  \ifthenelse{ \equal{#1}{1} }{i}{
  \ifthenelse{ \equal{#1}{2} }{ii}{
  \ifthenelse{ \equal{#1}{3} }{iii}{
  \ifthenelse{ \equal{#1}{4} }{iv}{
  \ifthenelse{ \equal{#1}{5} }{v}{
  \ifthenelse{ \equal{#1}{6} }{vi}{
  \ifthenelse{ \equal{#1}{7} }{vii}{
  \ifthenelse{ \equal{#1}{8} }{viii}{
  \ifthenelse{ \equal{#1}{9} }{ix}{
  \ifthenelse{ \equal{#1}{10} }{x}{
  \ifthenelse{ \equal{#1}{11} }{xi}{
  \ifthenelse{ \equal{#1}{12} }{xii}{
  \ifthenelse{ \equal{#1}{13} }{xiii}{
  \ifthenelse{ \equal{#1}{14} }{xiv}{
  \ifthenelse{ \equal{#1}{15} }{xv}{
  \ifthenelse{ \equal{#1}{16} }{xvi}{
  \ifthenelse{ \equal{#1}{17} }{xvii}{
  \ifthenelse{ \equal{#1}{18} }{xviii}{
    ERREUR,~MODIFIER~LA~MACRO~numtoi
  } } } } } } } } } } } } } } } } } }
}

\newcounter{cretraiti}
\newenvironment{ri}[1]{
\begin{list}{($\numtoi{\thecretraiti}$)}{
  \usecounter{cretraiti}
  \topsep=0.5ex
  \itemsep=0.3ex
  \labelsep=0.3em
  \parsep=0ex
  \listparindent=1em
  \settowidth{\labelwidth}{(#1)}
  \leftmargin=\labelwidth
}}{\end{list}
}

\newcommand{\Ascr}{\mathscr{A}}
\newcommand{\Bscr}{\mathscr{B}}
\newcommand{\Dscr}{\mathscr{D}}
\newcommand{\Escr}{\mathscr{E}}
\newcommand{\Fscr}{\mathscr{F}}
\newcommand{\Hscr}{\mathscr{H}}

\newcommand{\Oo}{\mathcal{O}}

\newcommand{\be}{\beta}

\newcommand{\accolades}[1]{\ensuremath{\left\{#1\right\}}}
\newcommand{\paa}[1]{\ensuremath{\accolades{#1}}}
\newcommand{\ins}[1]{{\ensuremath{\textrm{#1}}}}

\newcommand{\ra}{\rightarrow}
\newcommand{\id}{\ensuremath{\mathrm{id}}}

\begin{document}

\title{Construction fonctorielle de cat\'egorie de Frobenius}
\author{Beck Vincent}

\selectlanguage{francais}
\maketitle 

\begin{abstract}
\selectlanguage{francais}
Soient $\Ascr,\Bscr$ deux cat\'egories exactes telle que $\Ascr$ soit karoubienne et $M : \Bscr \ra \Ascr$ un foncteur exact.
Sous des hypoth\`eses d'adjonction pour $M$, on montre que les objets de $\Ascr$ qui sont facteurs directs d'objets de la forme $MY$ pour $Y \in \Bscr$
forment alors une cat\'egorie de Frobenius ce qui permet de d\'efinir par passage au quotient la cat\'egorie $M$-stable de $\Ascr$.
Par ailleurs, on sugg\`ere la construction d'une cat\'egorie $M$-stable pour $\Ascr,\Bscr$ des cat\'egories triangul\'ees et $M$ un foncteur triangul\'e.
On illustre cette derni\`ere notion par un th\'eor\`eme de Keller et Vossieck (voir~\cite{keller-vossieck}) qui relie les deux notions de cat\'egorie $M$-stable.
\vskip 0.5\baselineskip

\selectlanguage{english}
\begin{center}\noindent{\bf Abstract}
\noindent
\end{center}
{\bf Fonctorial Construction of Frobenius Categories. } Let $\Ascr,\Bscr$ be exact categories with $\Ascr$ karoubian and $M$ be an exact functor.
Under suitable adjonction hypotheses for $M$, we are able to show that the direct factors of the objects of $\Ascr$ of the form $MY$ with $Y \in \Bscr$
make up a Frobenius category which allow us to define an $M$-stable category for $\Ascr$ only by quotienting. In addition, we propose a
construction of an $M$-stable category for $\Ascr,\Bscr$ triangulated categories and $M$ a triangulated functor.
We illustrate this notion with a theorem of Keller and Vossieck (see~\cite{keller-vossieck}) which links the two notions of $M$-stable category.
\end{abstract}


\begin{msc} 18E10 - 18E30 \end{msc}

\selectlanguage{english}
\section*{Abridged English version}
In this article, we study the following situation : let us consider two additive categories $\Ascr,\Bscr$ and an additive functor $M : \Bscr \ra \Ascr$;
we are interested  in the quotient category of $\Ascr$ by the objects which are the direct factors of the objects $MY$ with $Y \in \Bscr$.
We first consider the case where $\Ascr$ is an exact category. In this case, we construct a structure of Frobenius category on $\Ascr$.
For this Frobenius structure, the injectives-projectives are nothing else but the direct factors of the $MY$ for $Y \in \Bscr$.
Thus, we obtain a triangulated structure on the quotient which is called the $M$-stable category associated to $\Ascr$.
On the second hand, we study the case of a triangulated category $\Ascr$ and also define a notion $M$-stable category in this case.
We finally link these two notions of $M$-stable category.

We begin with the case of an exact category $\Ascr$. We consider two exact categories $(\Ascr,\Escr)$ and $(\Bscr,\Fscr)$ and an exact
functor $M : \Bscr \ra \Ascr$. Under suitable adjonction hypotheses on $M$ (see~proposition~\ref{prop-frobenius-eng}),
we will construct a new exact structure on $\Ascr$ which is a Frobenius structure.
For this, we begin with the lemma~\ref{lemma-exact-structure} which allow us to construct easily new exact structure
by an inverse image process.

\vskip2ex
\begin{lemma}[New exact structure]\label{lemma-exact-structure}
Let $(\Ascr,\Escr)$, $(\Bscr,\Fscr)$ be two exact categories and $L : \Ascr \ra \Bscr$ be an exact functor.
We define the family  $\Escr_L^0 = \paa{s \in \Escr,\ Ls \ins{ is a split exact sequence in } \Bscr}\,.$
Then $(\Ascr,\Escr_L^0)$ is an exact category.
\end{lemma}

\vskip1ex
The proof of this lemma is straightforward by using the minimal axioms of Keller~\cite{keller}.
We can now state our main result on exact category.

\vskip2ex
\begin{proposition}[Frobenius category]\label{prop-frobenius-eng}
Let $(\Ascr,\Escr)$, $(\Bscr,\Fscr)$ be two exact categories and $M : \Bscr \ra \Ascr$ be an exact functor.
Let us assume that $M$ has a left adjoint $L$ and a right adjoint $R$. We also assume the following properties.

\begin{ri}{(iii)}
\item $\Ascr$ is karoubian;

\item $L$ and $R$ are exact functors between the exact categories $\Ascr$ and $\Bscr$;

\item $\Escr_L^0= \Escr_R^0:= \Escr'$;

%
%
\item for every $f:X \ra X'$ such that there exists $\be$ verifying $\be Lf = \id_{LX}$, the morphism $f$ is an inflation of $\Ascr$ (and so of $(\Ascr,\Escr')$).

\item for every $f:X \ra X'$ such that there exists $\be$ verifying $Rf \be = \id_{RX}$, the morphism $f$ is a deflation of $\Ascr$ (and so of $(\Ascr,\Escr')$).
\end{ri}
Then $(\Ascr,\Escr')$ is a Frobenius category whose projectives-injectives are the direct factors of the objets $MY$ for $Y \in \Bscr$.
(these projectives-injectives are called the $M$-split objects of $\Ascr$). The stable category of $(\Ascr,\Escr')$ is called the $M$-stable category of $\Ascr$.
\end{proposition}

\vskip1ex
This result can be seen as the mirror reflection of the result of~\cite{grime}.
These hypotheses are verified when $\Ascr = A$-Mod. is the category of $A$-modules with $A$ a Frobenius algebra over $k$,
$\Bscr$ is the category of vector spaces of $k$ and $M = A \otimes\,? : \Bscr \ra \Ascr$ is the induction functor.
The hypotheses are also verified when $\Ascr = \Oo G$-Mod. is the category of $\Oo G$-modules (where $\Oo$ is a commutative ring and $G$ finite group)
$\Bscr$ is the category of $\Oo H$-modules (where $H$ is a subgroup of $G$)  $M = \ins{Ind}_H^G : \Bscr \ra \Ascr$ is the induction functor.

\vskip1ex
Let us now study triangulated categories and for this, we begin with the definition of the $M$-stable category in this case.

\vskip2ex
\begin{definition}[The $M$-stable category of a triangulated category]\label{dfn-categ-mstable-triangulee-eng} Let $\Ascr$ and $\Bscr$ be two triangulated categories
and $M : \Bscr \ra \Ascr$ be a triangulated functor with left and right adjoint. The $M$-stable of $\Ascr$ is, by definition, the triangulated quotient
$\Ascr/\langle M\ins{-split\,} \rangle$ where $\langle M\ins{-split\,} \rangle$ is the thick subcategory of $\Ascr$ generated by
the $MY$ for $Y \in \Bscr$.
\end{definition}

\vskip1ex
\begin{example}[Derived category and $M$-stable category]\label{ex-eng-categ-mstable-triangulee}
We assume hypotheses of proposition~\ref{prop-frobenius-eng} and denote by $\Hscr^b(M)$ the bounded homotopy category of the $M$-split objects of $\Ascr$.
Following~\cite{keller-vossieck} or~\cite{keller-derived}, we have $M\ins{Stab}_\Ascr = \Dscr^{b}(\Ascr,\Escr')/\Hscr^b(M)$.

In addition, let us denote by $\Fscr_M^0$ the exact structure on $\Bscr$ obtained as inverse image of the structure of split sequences
of $\Ascr$ by $M$. The adjunctions $(L,M)$ et $(M,R)$ extend to triangulated adjunctions which are still denoted by $(L,M)$ et $(M,R)$ between the
triangulated categories $\Dscr^b(\Ascr,\Escr')$ and $\Dscr^b(\Bscr,\Fscr_M^0)$.
Using lemma~2.4 of~\cite{balmer-schlichting}, we obtain that, for this functor $M$, the thick subcategory of $\Dscr^b(\Ascr,\Escr')$
generated by the $MY$ for $Y \in \Dscr^b(\Bscr,\Fscr)$ is nothing else but $\Hscr^b(M)$.

Thus the $M$-stable category of the triangulated category $\Dscr^b(\Ascr,\Escr')$ is also the $M$-stable category of the exact category $\Ascr$.
\end{example}

\selectlanguage{francais}
\section{Introduction}\label{sec-intro}
Cet article propose une \'etude de la situation suivante : on consid\`ere une cat\'egorie additive $\Ascr$ munie d'un foncteur additif $M : \Bscr \ra \Ascr$.
On cherche \`a construire une structure sur la cat\'egorie quotient de $\Ascr$ par les facteurs directs d'objets de la forme $MY$.
Dans la section~\ref{sec-frobenius}, on s'int\'eresse au cas o\`u $\Ascr$ est une cat\'egorie exacte et on construit une structure de cat\'egorie de Frobenius
dont les projectifs-injectifs sont les objets facteurs directs d'objets de la forme $MY$ (c'est la proposition~\ref{prop-frobenius}).
On obtient ainsi une structure triangul\'ee sur le quotient (c'est la d\'efinition~\ref{dfn-categ-mstable}). Dans la section~\ref{sec-triangulee}, on
\'etudie le cas o\`u $\Ascr$ est une cat\'egorie triangul\'ee (c'est la d\'efinition~\ref{dfn-categ-mstable-triangulee})
et on relie cette construction triangul\'ee \`a celle de la section~\ref{sec-frobenius}.

\section{Construction de cat\'egorie de Frobenius}\label{sec-frobenius}

\subsection{Nouvelle structure exacte}\label{ssec-categ-exacte}

Dans le lemme~\ref{lem-structure-exacte} qui suit,
on construit sur une cat\'egorie exacte une nouvelle structure exacte (extraite de celle de d\'epart) par un proc\'ec\'e d'image r\'eciproque par un foncteur exact
depuis une cat\'egorie exacte disposant de deux structures exactes.
Le corollaire~\ref{cor-exact-scinde} met en exergue la situation de rel\`evement de la structure minimale donn\'ee par les suites exactes scind\'ees.
C'est cette derni\`ere structure qui sera utilis\'ee par la suite dans la sous-section~\ref{ssec-categ-frobenius} pour construire des cat\'egories de Frobenius.

\vskip2ex
\begin{lemme}[Nouvelle structure exacte]\label{lem-structure-exacte}
Soient $(\Ascr,\Escr)$ et $(\Bscr,\Fscr)$ deux cat\'egories exactes au sens de Quillen~\cite{quillen} et Keller~\cite{keller} et $L : \Ascr \ra \Bscr$ un foncteur exact.
On suppose que $\Fscr'$ est sous-famille de $\Fscr$ telle que $(\Bscr,\Fscr')$ est encore une cat\'egorie exacte.
On d\'efinit la famille $\Escr'=\paa{s \in \Escr,\ Ls \in \Fscr'}$. Le couple $(\Ascr,\Escr')$ est une cat\'egorie exacte.
\end{lemme}

\vskip1ex
\noindent {\bf Preuve.} On reprend les axiomes minimaux pr\'esent\'es dans l'appendice A de l'article~\cite{keller}.
Par d\'efinition, $\Escr'$ est stable par isomorphisme.

{\it Ex$0$.} Par hypoth\`ese, la suite $0 \rightarrowtail 0 \twoheadrightarrow 0$ est une suite exacte de $\Escr$ dont l'image par $L$ est dans $\Fscr'$.
Ainsi, elle est dans $\Escr'$ : $\id_0$ est une d\'eflation de $\Escr'$.

{\it Ex$1$.} Soient $f :X \ra X',g : X' \ra X''$ deux d\'eflations de $\Escr'$. Ce sont en particulier des d\'eflations de $\Escr$.
On a donc une suite exacte de $\Escr$ de la forme $Y \rightarrowtail X \twoheadrightarrow X''$.
Montrons que son image par $L$ est dans $\Fscr'$.
Comme $Lf$ et $Lg$ sont des d\'eflations de $\Fscr'$ et comme $(\Bscr,\Fscr')$ est exact, on a
une suite exacte de $\Fscr'$ de la forme $Z \rightarrowtail LX \twoheadrightarrow LX''$.
Comme $Z$ est le noyau de $LgLf$, on en d\'eduit que cette suite exacte est isomorphe \`a $LY \rightarrowtail LX \twoheadrightarrow LX''$.
Ainsi $gf$ est un d\'eflation de $\Escr'$.

{\it Ex$2$.} Soient $d: Y \ra Z$ une d\'eflation de $\Escr'$ et $f : Z' \ra Z$. Il existe une d\'eflation $d'$ de $\Escr$
et $f'$ tels que le diagramme
$$\xymatrix{Y' \ar^{d'}[r] \ar^{f'}[d]& Z' \ar^f[d] \\ Y \ar^d[r]& Z }$$
soit cart\'esien. Montrons que $d'$ est une d\'eflation de $\Escr'$.
Comme $L : (\Ascr,\Escr)\ra (\Bscr,\Fscr)$ est exacte, l'image par $L$ de ce diagramme est encore cart\'esienne
puisque $Y' \rightarrowtail Z' \oplus Y \twoheadrightarrow Z$ est une suite exacte de $\Escr$ (voir~\cite{keller}).
Par ailleurs, comme $Ld$ est une d\'eflation de $\Fscr'$, il existe une d\'eflation $d''$ de $\Fscr'$ et $f''$ tels que le diagramme
$$\xymatrix{Y'' \ar^{d''}[r] \ar^{f''}[d]& LZ' \ar^{Lf}[d] \\ LY \ar^{Ld}[r]& LZ }$$
soit cart\'esien. Par la propri\'et\'e universelle du produit fibr\'e, on en d\'eduit un diagramme commutatif
$$\xymatrix{Y'' \ar^{d''}[r] \ar^{g}[d]& LZ' \ar@{=}[d] \\ LY' \ar^{Ld'}[r]& LZ' }$$
o\`u $g$ est un isomorphisme. Comme $d''$ est une d\'eflation de $\Fscr'$ et $d'$ une d\'eflation de $\Escr$, on peut compl\'eter le diagramme pr\'ec\'edent en
$$\xymatrix{U \ar[d]\ar[r]& Y'' \ar^{d''}[r] \ar^{g}[d]& LZ' \ar@{=}[d] \\LX' \ar[r] & LY' \ar^{Ld'}[r]& LZ'}$$
o\`u les fl\`eches verticales sont des isomorphismes, $X' \rightarrowtail Y' \twoheadrightarrow Z'$ une suite exacte de $\Escr$ et
la premi\`ere ligne, une suite exacte de $\Fscr'$. Ainsi $d'$ est bien une d\'eflation de $\Escr'$.

{\it Ex$2^{op}$.} On l'obtient de la m\^eme fa\c{c}on que l'axiome pr\'ec\'edent. \myqed

\vskip2ex
En consid\'erant l'ensemble des suites scind\'ees de $\Bscr$, on obtient ainsi le r\'esultat suivant.
\vskip1ex

\begin{corollaire}[Rel\`evement des suites scind\'ees]\label{cor-exact-scinde}
Soient $(\Ascr,\Escr)$ et $(\Bscr,\Fscr)$ deux cat\'egories exactes et $L : \Ascr \ra \Bscr$ un foncteur exact.
On d\'efinit la famille $\Escr_L^0=\paa{s \in \Escr,\ Ls \ins{ scind\'ee }}$. Le couple $(\Ascr,\Escr_L^0)$ est une cat\'egorie exacte.
\end{corollaire}

\subsection{Cat\'egorie de Frobenius}\label{ssec-categ-frobenius}

Dans la proposition~\ref{prop-frobenius} qui suit, on donne un cadre g\'en\'eral permettant de construire une structure de cat\'egorie de Frobenius \`a partir
de la donn\'ee d'un triplet $(L,M,R)$ de foncteurs adjoints et exacts. De plus, pour cette structure de Frobenius, les injectifs-projectifs sont faciles \`a
d\'ecrire : ce sont les objets facteurs directs d'objets de la forme $MY$. On peut alors d\'efinir la cat\'egorie $M$-stable comme la cat\'egorie stable au sens
de Happel~\cite{happel} de la cat\'egorie de Frobenius ainsi obtenue : c'est la d\'efinition~\ref{dfn-categ-mstable}.
Enfin, dans l'exemple~\ref{ex-categ-abelienne}, on montre que les nombreuses hypoth\`eses de la proposition~\ref{prop-frobenius}
peuvent \^etre obtenues de fa\c{c}on \'el\'ementaire dans le cadre ab\'elien.

\vskip2ex
\begin{proposition}[Cat\'egorie de Frobenius]\label{prop-frobenius}
Soient $(\Ascr,\Escr)$ et $(\Bscr,\Fscr)$ deux cat\'egories exactes et un foncteur exact $M : \Bscr \ra \Ascr$
ayant un adjoint \`a droite $R$ et un adjoint \`a gauche $L$. On fait les hypoth\`eses suivantes.

\begin{ri}{(iii)}
\item $\Ascr$ est karoubienne;

\item $L$ et $R$ sont des foncteurs exacts entre les cat\'egories exactes $\Ascr$ et $\Bscr$;

\item $\Escr_L^0= \Escr_R^0:= \Escr'$;

%
%
\item pour tout $f:X \ra X'$ tel qu'il existe $\be$ v\'erifiant $\be Lf = \id_{LX}$, le morphisme $f$ est une inflation de $\Ascr$ (et donc de $(\Ascr,\Escr')$).

\item pour tout $f : X \ra X'$ tel qu'il existe $\be$ v\'erifiant $Rf \be = \id_{RX}$, le morphisme $f$ est une d\'eflation de $\Ascr$ (et donc de $(\Ascr,\Escr')$).
\end{ri}
Alors $(\Ascr,\Escr')$ est une cat\'egorie de Frobenius dont les projectifs-injectifs sont les facteurs directs d'objets de la forme $MY$ pour $Y \in \Bscr$
(ces projectifs-injectifs sont appel\'es les objets $M$-scind\'es de $\Ascr$).
\end{proposition}

\vskip1ex
\noindent {\bf Preuve.} D'apr\`es le corollaire~\ref{cor-exact-scinde}, $(\Ascr,\Escr')$ est une cat\'egorie exacte.
Par d\'efinition, les injectifs de $(\Ascr,\Escr')$ sont les objets $I$ de $\Ascr$ tels que
pour toute inflation $i:X \ra X'$ de $\Ascr$ telle qu'il existe $\be$ v\'erifiant $\be Li = \id_{LX}$ et pour tout $f : X \ra I$, il existe
$g : X' \ra I$ tel que $f=gi$. L'hypoth\`ese~$(iv)$ assure que les injectifs de $(\Ascr,\Escr')$ sont les objets v\'erifiant
pour tout $i:X \ra X'$ tel qu'il existe $\be$ v\'erifiant $\be Li = \id_{LX}$ et pour tout $f : X \ra I$, il existe
$g : X' \ra I$ tel que $f=gi$. Comme~$\Escr$ est karoubienne, les injectifs $(\Ascr,\Escr')$ sont les facteurs directs d'objets de la forme $MY$ pour $Y \in \Bscr$
(voir~\cite[th\'eor\`eme 6.8]{broue} et~\cite[proposition~3.24]{beck-these}).
En suivant le m\^eme raisonnement, on obtient que les projectifs de $(\Ascr,\Escr')$ sont aussi les facteurs directs d'objets de la forme $MY$.
Enfin, les propri\'et\'es des adjonctions assurent que, pour tout $X$, l'unit\'e $\eta_X : X \ra MLX$ de l'adjonction $(L,M)$ est telle qu'il existe
$\be L\eta_X=\id_{LX}$. Ainsi, par~$(iv)$, $(\Ascr,\Escr')$ a assez d'injectifs et, de m\^eme, assez de projectifs. \myqed

\vskip2ex
Ce r\'esultat est en un certain sens le sym\'etrique de celui de Grime~\cite{grime}. 
On d\'efinit maintenant la notion de cat\'egorie $M$-stable de la cat\'egorie $\Ascr$ v\'erifiant les hypoth\`eses de la proposition pr\'ec\'edente
simplement comme la cat\'egorie stable de la cat\'egorie de Frobenius construite.
\vskip1ex
\begin{definition}[Cat\'egorie M-stable]\label{dfn-categ-mstable} Sous les hypoth\`eses de la proposition pr\'ec\'edente, la cat\'egorie stable
(au sens de Happel~\cite{happel}) de la cat\'egorie exacte $(\Ascr,\Escr')$ est appel\'ee la cat\'egorie $M$-stable de $\Ascr$ et est not\'ee~$M\ins{stab}_\Ascr$.
\end{definition}

\vskip2ex
On donne \`a pr\'esent quelques exemples concrets de la situation d\'ecrite dans la proposition~\ref{prop-frobenius} dans le cas o\`u $\Ascr$ est une
cat\'egorie ab\'elienne.

\vskip1ex
\begin{exemple}[Le cas ab\'elien]\label{ex-categ-abelienne}
On suppose que $\Ascr$ et $\Bscr$ sont des cat\'egories ab\'eliennes munies de leur structure maximale de cat\'egorie exacte.
L'hypoth\`ese~$(i)$ est automatiquement v\'erifi\'ee. Par ailleurs, la fid\'elit\'e de $L$ et $R$ donne les hypoth\`eses $(iv)$~et~$(v)$.
On retrouve ainsi l'hypoth\`ese~4.14 de~\cite{beck-these}.

Par ailleurs, remarquons que les hypoth\`eses~$(ii)$ et~$(iii)$ sont toujours v\'erifi\'ees lorsque $L=R$.

Soient $k$ est un corps et $A$ une $k$-alg\`ebre de Frobenius ($\dim_k A < \infty$ et $A \simeq A^*$ en tant que $A$-module).
On consid\`ere $\Ascr= A$-$\ins{Mod}$ la cat\'egorie des $A$-modules, $\Bscr = k$-$\ins{Ev}$ la cat\'egorie des $k$-espaces vectoriels et
le foncteur d'induction $M = A \otimes\, ? : \Bscr \ra \Ascr$. On a alors $L=R$ est le foncteur de restriction qui est fid\`ele et exact.
Les injectifs-projectifs pour la structure exacte introduite ci-dessus ne sont rien d'autre que les $A$-modules projectifs au sens usuel.

Soient $\Oo$ un anneau commutatif unitaire, $G$ un groupe fini et $H$ un sous-groupe fini. On consid\`ere $\Ascr=\Oo G$-\ins{mod.} la
cat\'egorie des $G$-modules sur l'anneau $\Oo$, $\Bscr = \Oo H$-\ins{mod.} la cat\'egorie des $H$-modules sur $\Oo$ et $M = \Oo G \otimes_{\Oo H}\, ?$
le foncteur d'induction. On a alors $L=R$ qui est le foncteur de restriction. Il est bien fid\`ele et exact.
\end{exemple}

\section{Cat\'egorie M-stable d'une cat\'egorie triangul\'ee}\label{sec-triangulee}

Dans cette section, on sugg\`ere une construction d'une cat\'egorie $M$-stable d'une cat\'egorie triangul\'ee.
Dans la remarque~\ref{rem-schanuel-triangule}, on montre que, dans la cat\'egorie $M$-stable, le d\'ecalage se calcule \og \`a la Schanuel \fg.
Enfin, dans l'exemple~\ref{ex-categ-mstable-triangulee}, on relie, \`a l'aide d'un th\'eor\`eme de Keller et Vossieck,
les deux notions de cat\'egorie $M$-stable : celle d'une cat\'egorie exacte et celle d'une cat\'egorie triangul\'ee.

\vskip1ex
\begin{definition}[Cat\'egorie M-stable d'une cat\'egorie triangul\'ee]\label{dfn-categ-mstable-triangulee} On consid\`ere la situation suivante :
$\Ascr$ et $\Bscr$ sont deux cat\'egories triangul\'ees, $M : \Bscr \ra \Ascr$ un foncteur triangul\'e qui admet un adjoint
\`a droite $R$ et \`a gauche $L$, tous deux triangul\'es. On d\'efinit la cat\'egorie $M$-stable de $\Ascr$ comme la cat\'egorie triangul\'ee
quotient $\Ascr/\langle M\ins{-scind\'e\,} \rangle$ o\`u $\langle M\ins{-scind\'e\,} \rangle$ d\'esigne la sous-cat\'egorie \'epaisse de $\Ascr$ engendr\'ee
par les objets de la forme $MY$ pour $Y \in \Bscr$.
\end{definition}

\vskip1ex
\begin{remarque}[Lemme de Schanuel triangul\'e]\label{rem-schanuel-triangule}
Dans le cadre de la d\'efinition pr\'ec\'edente, le d\'ecalage peut se calculer de la fa\c{c}on suivante :
pour $X \in \Ascr$, on consid\`ere un couple $(P,i)$ o\`u $P$ est un facteur direct d'un objet de la forme $MY$ et $i : X \ra P$ (par exemple $P=MLX$ et $i = \eta_X$).
On compl\`ete alors $i : X \ra P$ en un triangle $X \ra P \ra \Omega X \ra X[1]$. Dans le quotient $P$ s'annule et ainsi $\Omega X$ s'identifie au d\'ecal\'e de $X$.
Si on impose en plus \`a $i$ de v\'erifier l'existence d'un $\be$ tel que $\be Li= \id_{LX}$ alors l'identification pr\'ec\'edente est fonctorielle
(voir~\cite[proposition 4.4]{beck-these}).
\end{remarque}

\vskip1ex
\begin{exemple}[Cat\'egorie d\'eriv\'ee et cat\'egorie $M$-stable]\label{ex-categ-mstable-triangulee}
On se place dans le cadre des hypoth\`eses de la proposition~\ref{prop-frobenius}.
On note $\Hscr^b(M)$ la cat\'egorie homotopique born\'ee de la sous-cat\'egorie pleine de $\Ascr$ form\'ee des objets $M$-scind\'es de $\Ascr$.
D'apr\`es~\cite{keller-vossieck} ou~\cite{keller-derived}, on a $M\ins{Stab}_\Ascr = \Dscr^{b}(\Ascr,\Escr')/\Hscr^b(M)$.

Par ailleurs, on note $\Fscr'$ la structure exacte sur $\Bscr$ obtenue comme image r\'eciproque de la structure des suites exactes scind\'ees
de $\Ascr$ par le foncteur $M$ (voir le corollaire~\ref{cor-exact-scinde}). Les adjonctions $(L,M)$ et $(M,R)$ s'\'etendent en des adjonctions (triangul\'ees)
encore not\'ee $(L,M)$ et $(M,R)$ entre les cat\'egories triangul\'ees $\Dscr^b(\Ascr,\Escr')$ et $\Dscr^b(\Bscr,\Fscr')$.
Pour ce foncteur $M$, la sous-cat\'egorie \'epaisse de $\Dscr^b(\Ascr,\Escr')$
engendr\'ee par les objets de la forme $MY$ pour $Y \in \Dscr^b(\Bscr,\Fscr)$ (qu'on note $\langle M\ins{-scind\'e\,} \rangle$), est en fait $\Hscr^b(M)$.
En effet, comme $\Ascr$ est karoubienne, la sous-cat\'egorie pleine de $\Ascr$ form\'ee des objets $M$-scind\'es l'est aussi. Le lemme~2.4 de~\cite{balmer-schlichting}
assure que $\Hscr^b(M)$ l'est aussi. Par ailleurs, la cat\'egorie $\langle M\ins{-scind\'e\,} \rangle$ contient bien entendu les complexes
ayant une seule composante non nulle qui est $M$-scind\'ee et donc, par extension, elle contient aussi $\Hscr^b(M)$.

Ainsi la cat\'egorie $M$-stable de la cat\'egorie triangul\'ee $\Dscr^b(\Ascr,\Escr')$ co\"{\i}ncide  avec la cat\'egorie $M$-stable de la cat\'egorie $\Ascr$.
\end{exemple}

\section*{Remerciements}
Je remercie grandement Bernhard Keller pour ses conseils avis\'es ainsi que pour les nombreuses r\'ef\'erences bibliographiques qu'il m'a donn\'ees.

\end{document}